\documentclass[11pt, twoside, leqno]{article}

\usepackage{amsmath,amsthm}
\usepackage{amssymb}

\usepackage{graphicx}

\usepackage[T1]{fontenc}

\usepackage[margin=1in]{geometry}
\usepackage{latexsym, amsfonts, amsmath,tikz}
\usepgflibrary{arrows}
\newcommand{\be}{\begin{equation}}
\newcommand{\ee}{\end{equation}}
\newcommand{\bea}{\begin{eqnarray}}
\newcommand{\eea}{\end{eqnarray}}
\newcommand{\bean}{\begin{eqnarray*}} 
\newcommand{\eean}{\end{eqnarray*}}

\newcommand{\brray}{\begin{array}}
\newcommand{\erray}{\end{array}}
\newcommand{\ben}{\begin{equation}{nonumber}}
\newcommand{\een}{\end{equation}{nonumber}}

\newtheorem{dfn}{Definition}[section]
\newtheorem{thm}[dfn]{Theorem}
\newtheorem{lmma}[dfn]{Lemma}
\newtheorem{ppsn}[dfn]{Proposition}
\newtheorem{crlre}[dfn]{Corollary}
\newtheorem{xmpl}[dfn]{Example}
\newtheorem{rmrk}[dfn]{Remark}
\newcommand{\bdfn}{\begin{dfn}}
\newcommand{\bthm}{\begin{thm}}
\newcommand{\blmma}{\begin{lmma}}
\newcommand{\bppsn}{\begin{ppsn}}
\newcommand{\bcrlre}{\begin{crlre}}
\newcommand{\bxmpl}{\begin{xmpl}}
\newcommand{\brmrk}{\begin{rmrk}}
\newcommand{\edfn}{\end{dfn}}
\newcommand{\ethm}{\end{thm}}
\newcommand{\elmma}{\end{lmma}}
\newcommand{\eppsn}{\end{ppsn}}
\newcommand{\ecrlre}{\end{crlre}}
\newcommand{\exmpl}{\end{xmpl}}
\newcommand{\ermrk}{\end{rmrk}}

\newcommand{\cla}{{\cal A}}
\newcommand{\clb}{{\cal B}}
\newcommand{\clc}{{\cal C}}
\newcommand{\cld}{{\cal D}}
\newcommand{\cle}{{\cal E}}
\newcommand{\clf}{{\cal F}}

\newcommand{\clh}{{\cal H}}
\newcommand{\cli}{{\cal I}}
\newcommand{\clj}{{\cal J}}

\newcommand{\cln}{{\cal N}}
\newcommand{\clo}{{\cal O}}

\newcommand{\clq}{{\cal Q}}

\newcommand{\clt}{{\cal T}}

\newcommand{\clv}{{\cal V}}
\newcommand{\clw}{{\cal W}}

\newcommand{\clz}{{\cal Z}}

\def\a*{{\cal A}_{h,*}}
\def\B{{\cal B}(h)}
\def\B1{{\cal B}_1(h)}
\def\b{{\cal B}^{\rm s.a.}(h)}
\def\b1{{\cal B}^{\rm s.a.}_1(h)}

\newcommand{\ot}{\otimes}

\newcommand{\raro}{\rightarrow}

\def \qed {$\Box$}

\newcommand{\midarrow}{\tikz \draw[-triangle 90] (0,0) -- +(.1,0);}

\def\a*{{\cal A}_{h,*}}
\def\B{{\cal B}(h)}
\def\B1{{\cal B}_1(h)}
\def\b{{\cal B}^{\rm s.a.}(h)}
\def\b1{{\cal B}^{\rm s.a.}_1(h)}

\begin{document}

	\baselineskip=17pt
	

	\title{Scalar curvature of a Levi-Civita connection on Cuntz algebra with three generators}
	
	\author{Soumalya Joardar\\
	Jawaharlal Nehru Center for Advanced Scientific Research\\ 
		Jakkur, Bangalore-560064\\
		Karnataka\\
		India\\
		E-mail: soumalya.j@gmail.com}	
	
	\date{}
	
	\maketitle
	\renewcommand{\thefootnote}{}
	
	
	\renewcommand{\thefootnote}{}
	\begin{abstract}
	A differential calculus on Cuntz algebra with three generators coming from the action of rotation group in three dimensions is introduced. The differential calculus is shown to satisfy Assumptions {\bf I-IV} of \cite{levi} so that Levi-Civita connection exists uniquely for any pseudo-Riemannian metric in the sense of \cite{levi}. Scalar curvature is computed for the Levi-Civita connection corresponding to the canonical bilinear metric.
	\end{abstract}
	
	
	\section{Introduction}
	Notion of connection and curvature is important in any form of geometry, be it classical or noncommutative. In noncommutative geometry, in last few years varied notions of curvature have been introduced. Broadly speaking there seems to be two main avenues to define curvature of a noncommutative space. One is to define directly Ricci and scalar curvature via some asymptotic expansion of noncommutative Laplacian (see for example \cite{moscovici}, \cite{Sitarz}, \cite{khalkhali}) and the other is to prove existence and uniqueness of Levi-Civita connection on the module of one forms and subsequently compute the curvature operator and scalar curvatures (see \cite{Frolich}, \cite{rosenberg}, \cite{Landi2}). The last one is more algebraic and has the advantage of finding the curvature operator directly. But the difficulty with this approach is to prove the existence and uniqueness of Levi-Civita connection on noncommutative spaces. Recently in \cite{levi}, the authors gave some sufficient conditions on the differential calculus on noncommutative spaces to ensure existence and uniqueness of Levi-Civita connection with respect to any pseudo-Riemannian metric in some sense. In a follow up paper (\cite{levi2}) the authors have given examples of noncommutative spaces admitting unique Levi-Civita connection in the sense of \cite{levi}.\\
	\indent In the context of noncommutative geometry it seems that any formulation of Levi-Civita connection invites automatically examples of noncommutative spaces accommodating such formulation leading to the existence and uniqueness theorem of Levi-Civita connection. A formulation also has to bode well in terms of computability of Ricci and scalar curvature. In this paper we introduce a differential calculus on Cuntz algebra with three generators coming from the natural action of the rotation group in three dimensions. Then we show that the differential calculus admits a unique Levi-Civita connection in the sense of \cite{levi}. We compute the Ricci curvature and scalar curvature along the lines of \cite{levi} with respect to the canonical metric. The scalar curvature turns out to be a negative constant times identity. It is worth mentioning that noncommutative geometry of Cuntz algebra has been studied in more analytic set up (see \cite{Rennie1} for example). Absence of any faithful trace on Cuntz algebra poses considerable amount of difficulty there. This paper attempts in a sense more algebraic study of the geometry of Cuntz algebra. 
	\section{Preliminaries}
	\subsection{Existence and uniqueness of Levi-Civita connection on a class of modules of one-forms}
	\label{prel}
	We begin by recalling the definition of a differential calculus on a $\ast$-algebra $\cla$ over $\mathbb{C}$. 
	\bdfn
	\label{differential}
	A differential calculus on a $\ast$-algebra $\cla$ is a pair $(\Omega(\cla),d)$ such that\\
	(i) $\Omega(\cla)$ is an $\cla-\cla$-bimodule.\\
	(ii) $\Omega(\cla)=\oplus_{i\geq 0}\Omega^{i}(\cla)$, where $\Omega^{0}(\cla)=\cla$ and each $\Omega^{i}(\cla)$ is an $\cla-\cla$-bimodule.\\
	(iii) There is a bimodule map $m:\Omega^{i}(\cla)\ot_{\cla}\Omega^{j}(\cla)\subset \Omega^{i+j}(\cla)$ for all $i,j$.\\
	(iv) $d:\Omega^{i}(\cla)\raro\Omega^{i+1}(\cla)$ satisfies $d(\omega.\eta)=d\omega.\eta+(-1)^{\rm deg \omega}\omega.d\eta$ and $d^2=0$.\\
	(v) $\Omega^{i}(\cla)$ is spanned by $da_{0}...da_{i}a_{i+1}$.
	\edfn 
	Let $(\Omega(\cla),d)$ be such a differential calculus on an algebra $\cla$.
	\bdfn
	A $\mathbb{C}$-linear map $\nabla:\Omega^{1}(\cla)\raro\Omega^{1}(\cla)\ot_{\cla}\Omega^{1}(\cla)$ is said to be a (right) connection on the $\cla-\cla$-bimodule $\Omega^{1}(\cla)$ if
	\begin{displaymath}
	\nabla(\omega.a)=\nabla(\omega)a+\omega\ot da, \omega\in\Omega^{1}(\cla), a\in\cla.
	\end{displaymath}
	\edfn
	\bdfn
	\label{torsionless}
	A connection $\nabla$ is said to be torsionless if the right linear map $T_{\nabla}=m\circ\nabla+d$ is equal to zero.
	\edfn
	For a bimodule $\cle$ over $\cla$, $\clz(\cle)$ will denote the center of the module i.e.
	\begin{displaymath}
	\clz(\cle)=\{e\in\cle:a.e=e.a \ for \ all \ a\in\cla\}.
	\end{displaymath}
	Similarly $\clz(\cla)$ will denote the center of an algebra $\cla$.
	Now we are going to state a few assumptions on a given differential calculus over a $\ast$-algebra $\cla$ following \cite{levi} which enables one to prove existence and uniqueness of Levi-Civita connection on the space of one forms as in \cite{levi}. In the following we shall denote the bimodule of one forms by $\cle$.\vspace{0.1in}\\
	{\bf Assumption I}: The $\cla-\cla$ bimodule $\cle$ is finitely generated projective right $\cla$-module. Moreover the map $u^{\cle}:\clz(\cle)\ot_{\clz(\cla)}\cla\raro\cle$ given by $u^{\cle}(\sum_{i}e_{i}\ot a_{i})=\sum e_{i}a_{i}$ is an isomorphism of vector spaces.\vspace{0.1in}\\
	{\bf Assumption II}: The right $\cla$-module $\cle\ot_{\cla}\cle$ admits a splitting $\cle\ot_{\cla}\cle={\rm ker}(m)\oplus \clf$ where $\clf\cong {\rm Im}(m)$.\vspace{0.1in}\\
	{\bf Assumption III}: If we denote the idempotent projecting onto ${\rm ker}(m)$ by $P_{\rm sym}$, then the map $\sigma:=(2P_{\rm sym}-1)$ is such that $\sigma(\omega\ot\eta)=\eta\ot\omega$ for $\omega,\eta\in\clz(\cle)$.\vspace{0.1in}\\
	Before stating the ${\bf IV}$th and final assumption, we recall the definition of a pseudo-Riemannian metric on a differential calculus over a $\ast$-algebra $\cla$ which satisfies Assumptions {\bf I-III}.
	\bdfn
	\label{pseudometric}
	A pseudo-Riemannian metric $g$ on $\cle$ is an element of ${\rm Hom}_{\cla}(\cle\ot_{\cla}\cle,\cla)$ such that $g\circ\sigma=g$ and it is non degenarate in the sense that the map $V_{g}:\cle\raro\cle^{\ast}$ defined by $V_{g}(\omega)(\eta)=g(\omega\ot\eta)$ is an isomorphism of right $\cla$-modules. 
	\edfn 
	Now we are ready to state the final assumption.\\
		{\bf Assumption IV}: There is a {\bf bilinear} pseudo-Riemannian metric $g\in {\rm Hom}_{\cla}(\cle\ot_{\cla}\cle,\cla)$. 
	\bdfn (see Proposition 2.12 qnd definition 2.13 of \cite{levi})
	Suppose for a $\ast$-algebra $\cla$, the differential calculus satisfies assumptions ${\bf I-IV}$. Then a connection $\nabla$ on the space of one forms is said to be unitary with respect to a pseudo-Riemannian metric $g$ if $dg=\Pi_{g}(\nabla)$, where $\Pi_{g}(\nabla):\cle\ot_{\cla}\cle\raro \cle$ is defined as
	\begin{displaymath}
	\Pi_{g}(\nabla)(\omega\ot\eta)=(g\ot{\rm id})\sigma_{23}(\nabla(\omega)\ot\eta+\nabla(\eta)\ot\omega).
	\end{displaymath}
\edfn
Recall the definition of the right $\cla$-linear map $\Phi_{g}:{\rm Hom}(\cle,\cle\ot^{\rm sym}_{\cla}\cle)\raro{\rm Hom(\cle\ot^{\rm sym}_{\cla}\cle,\cle)}$ (see Theorem 2.14 of \cite{levi}). We state the following theorem from \cite{levi} (see theorem 3.12 and 2.14):
\bthm
\label{gen}
For a $\ast$-algebra $\cla$ such that the differential calculus over $\cla$ satisfies assumptions ${\bf I-IV}$, a unique Levi-Civita connection (a connection which is torsionless and unitary) exists. Moreover for any torsionless connection $\nabla_{0}$, the Levi-Civita connection is given by
\begin{displaymath}
\nabla=\nabla_{0}+L,
\end{displaymath}
where $L=\Phi_{g}^{-1}(dg-\Pi_{g}(\nabla_{0}))$.
\ethm
Let us recall some definitions which will be used later in the paper. For a right $\cla$ module $\cle$, $\cle^{\ast}$ will denote the right $\cla$-module ${\rm Hom}_{\cla}(\cle,\cla)$. For two bimodules $\cle,\clf$, the map $\zeta_{\cle,\clf}:\cle\ot_{\cla}\clf^{\ast}\raro {\rm Hom}_{\cla}(\clf,\cle)$ is defined as
\begin{eqnarray}
\label{zeta}
\zeta_{\cle,\clf}(\sum_{i}e_{i}\ot\phi_{i})(f)=\sum_{i}e_{i}(\phi_{i}(f)).
\end{eqnarray}
When both $\cle$ and $\clf$ are finitely generated, projective, $\zeta_{\cle,\clf}$ is an {\bf isomorphism}. \\
Let $\cle$ be an $\cla-\cla$ bimodule satisfying assumptions $I-IV$. Then for any $\cla-\cla$ bimodule $\clf$, the map $T^{L}_{\cle,\clf}:\cle\ot_{\cla}\clf\raro\clz(\cle)\ot_{\clz(\cla)}\clf$ given by
\begin{eqnarray}
\label{TL}
T^{L}_{\cle,\clf}:=((u^{\cle})^{-1}\ot{\rm id}_{\clf}),
\end{eqnarray}	
	is a right $\cla$-linear isomorphism. Similarly we have the map $T^{R}_{\cle,\clf}$ (see Proposition 4.3) which is a left $\cla$-linear isomorphism between $\clf\ot_{\cla}\cle$ and $\clf\ot_{\clz(\cla)}\clz(\cle)$. Now recall the map $\rho:\cle\ot_{\cla}\cle^{\ast}\raro\cle^{\ast}\ot_{\cla}\cle$ defined as $\rho:=(T^{R}_{\cle,\cle^{\ast}})^{-1}{\rm flip}T^{L}_{\cle,\cle^{\ast}}$, where ${\rm flip}:\clz(\cle)\ot_{\clz(\cla)}\cle^{\ast}\raro \cle^{\ast}\ot_{\clz(\cla)}\clz(\cle)$ is the map ${\rm flip}(e^{\prime}\ot\phi)=\phi\ot e^{\prime}$. The map ${\rm ev}:\cle^{\ast}\ot\cle\raro\cla$ is given by ${\rm ev}(\sum_{i}\phi_{i}\ot e_{i})=\sum_{i}\phi_{i}(e_{i})$.\\
	\indent We also recall the map $H:\cle\ot_{\cla}\cle\raro\cle\ot_{\cla}\cle\ot_{\cla}\cle$ (see Lemma 4.1 of \cite{levi}) for a connection $\nabla$.
	\begin{eqnarray}
	\label{H}
	H(\omega\ot\eta)=(1-P_{\rm sym})_{23}(\nabla(\omega)\ot\eta)+\omega\ot Q^{-1}(d\eta),
	\end{eqnarray}
	where $Q:{\rm Im}(1-P_{\rm sym})\raro {\rm Im}(m)$ is the isomorphism from assumption II.
	\section{Connes' space of forms for Cuntz algebra with three generators} 
	Let $\cla$ be a $\ast$-algebra. We call a triple $(\cla,\clh,\cld)$ where $\cla$ is faithfully represented on $\clb(\clh)$ and $\cld$ is apriori an unbounded operator, a Dirac triple if $[\cld,a]\in\clb(\clh)$ for all $a\in\cla$. Note that a Dirac triple in our sense is called spectral triple in general. But a spectral triple comes with some kind of summability or compactness criterion. We don't need this in the paper and that is why we are content with a Dirac triple. For a $\ast$-algebra $\cla$, recall the reduced universal differential algebra $(\Omega^{\bullet}(\cla):=\oplus_{k} \Omega^{k}(\cla),\delta)$ from \cite{Connes}. Given a Dirac triple $(\cla,\clh,\cld)$ over $\cla$, there is a well defined $\ast$-representation $\Pi$ of $\Omega^{\bullet}(\cla)$ on $\clb(\clh)$ given by (see \cite{Landi})
	\begin{displaymath}
	\Pi(a_{0}\delta a_{1}...\delta a_{k})=a_{0}[\cld,a_{1}]...[\cld,a_{k}], a_{0},...,a_{k}\in\cla.
	\end{displaymath}
	Let $J^{k}_{0}=\{\omega\in\Omega^{k}(\cla):\Pi(\omega)=0\}$. The Connes' space of $k$-forms is defined to be 
	\begin{displaymath}
	\Omega^{k}_{\cld}(\cla)=\Pi(\Omega^{k}(\cla))/\Pi(\delta J^{k-1}_{0}).
	\end{displaymath}
	$\Pi(\delta J^{k-1}_{0})$ is a two sided ideal of $\Pi(\Omega^{k}(\cla))$ and is called the space of junk forms. For any element $\omega\in\Omega^{k}(\cla)$, if we denote the image of $\Pi(\omega)$ in $\Omega^{k}_{\cld}(\cla)$ by $\overline{\Pi(\omega)}$, then it can be shown that $(\Omega^{\bullet}_{\cld}(\cla):=\oplus_{k}\Omega^{k}_{\cld}(\cla),d)$ where $d$ is defined as $d\overline{\Pi(\omega)}:=\overline{\Pi(\delta\omega)}$, satisfies the conditions of definition \ref{differential}. 
\subsection{A Dirac triple on $\clo_{3}$: Connes' space of forms}
We begin this subsection by recalling the definition of Cuntz algebra with $n$-generators. The Cuntz algebra $\clo_{n}$ with $n$-generators is defined to be the universal $C^{\ast}$-algebra generated by $n$ elements $\{S_{i}\}_{i=1,\ldots,n}$ satisfying
\begin{eqnarray*}
&& S_{i}^{\ast}S_{i}=1, i=1,\ldots,n\\
&& \sum_{j=1}^{n}S_{j}S_{j}^{\ast}=1.
\end{eqnarray*} We shall limit ourselves to the Cuntz algebra $\clo_{3}$ with three generators in this paper. For $A=((a_{ij}))_{i,j=1,2,3}\in SO(3)$, define $\alpha_{A}(S_{i})=\sum_{j=1}^{3}a_{ij}S_{j}$ for $i=1,2,3$. Then 
	\begin{eqnarray*}
		\alpha_{A}(S_{i})^{\ast}\alpha_{A}(S_{i})&=&\sum_{j,k}a_{ij}a_{ik}S_{j}^{\ast}S_{k}\\
		&=& \sum_{j}a_{ij}a_{ij}\\
		&=& 1.
	\end{eqnarray*}
	Similarly we have $\sum_{i=1}^{3}\alpha_{A}(S_{i})\alpha_{A}(S_{i})^{\ast}=1$ so that for a fixed $A\in SO(3)$, we have a well defined $C^{\ast}$-homomorphism $\alpha_{A}:\clo_{3}\raro\clo_{3}$ which in turn gives a $C^{\ast}$-dynamical system $(\clo_{3},\alpha,SO(3))$. Now let us describe a Dirac triple arising from the above $C^{\ast}$-dynamical system. To that end recall that $SO(3)$ has three dimensional Lie-algebra $so(3)$. Let us choose three basis elements
		\begin{displaymath}
		X_{1}=\begin{bmatrix} 
		0 & 0 & 0\\           
		0 & 0 & -1\\           
		0 & 1 & 0             
		\end{bmatrix},
		X_{2}=\begin{bmatrix} 
		0 & 0 & -1\\           
		0 & 0 & 0\\           
		1 & 0 & 0             
		\end{bmatrix},
		X_{3}=\begin{bmatrix} 
		0 & 1 & 0\\           
		-1 & 0 & 0\\           
		0 & 0 & 0             
		\end{bmatrix},
		\end{displaymath} which defines three derivations on $\clo_{3}$ as
		\begin{displaymath}
		\partial_{i}(a)=\frac{d}{d\theta}\alpha_{{\rm exp}(\theta X_{i})}(a)|_{\theta=0}, a\in\clo_{3}, i=1,2,3.
		\end{displaymath}
		Note the following:
		\begin{displaymath}
		{\rm exp}(\theta X_{1})=\begin{bmatrix}
		1 & 0 & 0\\
		0 & cos(\theta) & -sin(\theta)\\
		0 & sin(\theta) & cos(\theta)
		\end{bmatrix},
		{\rm exp}(\theta X_{2})=\begin{bmatrix}
		cos(\theta) & 0 & -sin(\theta)\\
		0 & 1 & 0\\
		sin(\theta) & 0 & cos(\theta)
		\end{bmatrix},
		\end{displaymath}
		\begin{displaymath}
		{\rm exp}(\theta X_{3})=\begin{bmatrix}
		cos(\theta) & sin(\theta) & 0\\
		-sin(\theta) & cos(\theta) & 0\\
		0 & 0 & 1
		\end{bmatrix}.
		\end{displaymath}
		Then it is easy to see the actions of the derivations on the generators of $\clo_{3}$. They are given by: 
	\begin{eqnarray*}
		&&\partial_{1}(S_{1})=0, \partial_{1}(S_{2})=-S_{3},\partial_{1}(S_{3})=S_{2}\\
		&&\partial_{2}(S_{1})=-S_{3},\partial_{2}(S_{2})=0,\partial_{2}(S_{3})=S_{1}\\
		&&\partial_{3}(S_{1})=S_{2},\partial_{3}(S_{2})=-S_{1},\partial_{3}(S_{3})=0.
	\end{eqnarray*}
	It is clear that the derivations are actually $\ast$-derivations. The commutation relations are given by
	\begin{eqnarray}
	\label{commutations}
	[\partial_{1},\partial_{2}]=-\partial_{3},[\partial_{2},\partial_{3}]=\partial_{1},[\partial_{1},\partial_{3}]=-\partial_{2}.
	\end{eqnarray}\\
	\indent We denote the unique faithful KMS state of $\clo_{3}$ by $\tau$. Let $\clh$ be the Hilbert space $L^{2}(\clo_{3},\tau)\ot\mathbb{C}^{3}$. Then we define a Dirac triple $(\clo_{3},\clh,\cld)$ where $\pi:\clo_{3}\raro\clb(\clh)$ is given by $\pi(a):=a\ot \mathbb{I}$ and $\cld$ is the densely defined operator $\sum_{i=1}^{3}\partial_{i}\ot\sigma_{i}$, $\sigma_{i}$'s being $3\times 3$ Pauli Spin matrices satisfying $\sigma_{i}^{2}=\mathbb{I}$ and $\sigma_{i}\sigma_{j}=-\sigma_{j}\sigma_{i}$ for $i\neq j$. Using the derivation properties of $\partial_{i}$'s, note that for all $a\in\clo_{3}$,
	\begin{displaymath}[\cld,\pi(a)]=\sum_{i=1}^{3}(\partial_{i}(a))\ot\sigma_{i}.\end{displaymath}
	Hence $(\clo_{3},\clh,\cld)$ is in deed a Dirac triple.
	\bppsn
	\label{1forms}
	Connes' space of one forms $\Omega^{1}_{\cld}(\clo_{3})$ is a free module of rank $3$.
	
	\eppsn 
	{\it Proof}:\\
	Recall the differential $\delta$ on the reduced universal differential algebra $\Omega^{\bullet}(\clo_{3})$ as well as the well defined representation $\Pi:\Omega^{\bullet}(\clo_{3})\raro\clb(\clh)$. We have
	\begin{displaymath}
	\Omega^{1}_{\cld}(\clo_{3})=\{\Pi(\sum_{i} a_{i}\delta b_{i})\subset\clb(\clh):a_{i},b_{i}\in\clo_{3}\}.
	\end{displaymath}
 For $a_{i},b_{i}\in\clo_{3}$, 
	\begin{eqnarray*}\Pi(\sum_{i} a_{i}\delta b_{i})&=&\sum_{i}(a_{i}\ot\mathbb{I})[\cld,\pi(b_{i})]\\
	&=& \sum_{j=1}^{3}(\sum_{i}a_{i}\partial_{j}(b_{i})\ot \sigma_{j}),
	\end{eqnarray*}
	which proves that $\Omega^{1}_{\cld}(\clo_{3})\subset\clo_{3}\oplus\clo_{3}\oplus\clo_{3}$. We shall show that $\Omega^{1}_{\cld}(\clo_{3})=\clo_{3}\oplus\clo_{3}\oplus\clo_{3}$. To that end observe that 
	\begin{eqnarray*}\Pi(S_{1}^{\ast}\delta S_{2})&=&(S_{1}^{\ast}\ot\mathbb{I})[\cld,S_{2}]\\
	&=& -S_{1}^{\ast}S_{1}\ot\sigma_{3}\\
	&=& -1\ot\sigma_{3}.
	\end{eqnarray*}
	Similarly it can be shown that $\Pi(S_{1}^{\ast}\delta S_{3})=1\ot \sigma_{2}$ and $\Pi(S_{2}^{\ast}\delta S_{3})=1\ot\sigma_{1}$. Hence
	\begin{displaymath}
	\Omega^{1}_{\cld}(\clo_{3})=\clo_{3}\oplus\clo_{3}\oplus\clo_{3}.
	\end{displaymath}
	So $\Omega_{\cld}^{1}(\clo_{3})$ is a free module of rank $3$. Note that the freeness follows from the linear independence of $\sigma_{1},\sigma_{2}, \sigma_{3}$.\qed
	\bppsn
	\label{2forms}
	$\Omega^{2}_{\cld}(\clo_{3})$ is also a free module of rank $3$.
	\eppsn
	{\it Proof}:\\
	By definition $\Omega^{2}_{\cld}(\clo_{3})=\Pi(\Omega^{2}(\clo_{3}))/\Pi(\delta J_{0}^{1})$. We claim that 
	\begin{displaymath}
	\Pi(\Omega^{2}(\clo_{3}))=\clo_{3}\oplus\clo_{3}\oplus\clo_{3}\oplus\clo_{3}.
	\end{displaymath}
	To prove the above, let $\sum_{i}a_{i}\delta b_{i}\delta c_{i}\in\Omega^{2}(\clo_{3})$. Then
	\begin{eqnarray*}
	\Pi(\sum_{i}a_{i}\delta b_{i}\delta c_{i})&=& \sum_{i}(a_{i}\ot\mathbb{I})[\cld,b_{i}][\cld,c_{i}]\\
	&=& \sum_{i}(a_{i}\ot\mathbb{I})(\sum_{j=1}^{3}\partial_{j}(b_{i})\ot\sigma_{j})(\sum_{k=1}^{3}\partial_{k}(c_{i})\ot\sigma_{k})\\
	&=&\sum_{i} (a_{i}\sum_{j=1}^{3}\partial_{j}(b_{i})\partial_{j}(c_{i})\ot\mathbb{I}+\sum_{j<k}a_{i}(\partial_{j}(b_{i})\partial_{k}(c_{i})-\partial_{k}(b_{i})\partial_{j}(c_{i}))\ot\sigma_{j}\sigma_{k}),
	\end{eqnarray*}
	which proves $\Pi(\Omega^{2}(\clo_{3}))\subset\clo_{3}\oplus\clo_{3}\oplus\clo_{3}\oplus\clo_{3}$. To prove the equality note that $\Pi(\delta(S_{1}^{\ast})\delta(S_{1}))=\sum_{j=1}^{3}\partial_{j}(S_{1}^{\ast})\partial_{j}(S_{1})\ot\mathbb{I}=2\ot\mathbb{I}$. Also,
	\begin{displaymath}
	\Pi(\delta S_{2}^{\ast}\delta S_{1})=\sum_{j<k}(\partial_{j}(S_{2}^{\ast})\partial_{k}(S_{1})-\partial_{k}(S_{2}^{\ast})\partial_{j}(S_{1}))=1\ot\sigma_{1}\sigma_{2}.
	\end{displaymath}
	Similarly choosing appropriate elements from $\clo_{3}$, we can establish the equality. We are left with finding the junk forms. First, we shall show that $\Pi(\delta J_{0}^{1})\subset \clo_{3}\ot\mathbb{I}$. To that end choose $\omega=\sum_{i}a_{i}\delta b_{i}$ so that $\Pi(\omega)=0$. $\Pi(\omega)=0$ implies the following:
		\begin{eqnarray*}
		&&\sum_{i}(a_{i}\ot\mathbb{I})(\sum_{j=1}^{3}\partial_{j}(b_{i})\ot\sigma_{j})=0\\
		&\Rightarrow& \sum_{j}(\sum_{i}a_{i}\partial_{j}(b_{i}))\ot{\sigma_{j}}=0.
		\end{eqnarray*}
		Using the linear independence of $\sigma_{j}'s$, we get $\sum_{i}a_{i}\partial_{j}(b_{i})=0$ for $j=1,2,3$. Note that $\delta\omega=\sum_{i}\delta a_{i}\delta b_{i}$. Hence 
		\begin{displaymath}
		\Pi(\delta\omega)=\sum_{i}(\sum_{j=1}^{3}\partial_{j}(a_{i})\partial_{j}(b_{i})\ot\mathbb{I}+\sum_{j<k}(\partial_{j}(a_{i})\partial_{k}(b_{i})-\partial_{k}(a_{i})\partial_{j}(b_{i}))\ot\sigma_{j}\sigma_{k}).
		\end{displaymath}
		$\sum_{i}a_{i}\partial_{j}(b_{i})=0$ implies that for all $j,k=1,2,3$, $\sum_{i}\partial_{k}(a_{i})\partial_{j}(b_{i})=-a_{i}\partial_{k}\partial_{j}(b_{i})$. Hence
		\begin{displaymath}
		\Pi(\delta\omega)= \sum_{i}(\sum_{j=1}^{3}\partial_{j}(a_{i})\partial_{j}(b_{i})\ot\mathbb{I}-\sum_{j<k}(a_{i}[\partial_{j},\partial_{k}](b_{i})\ot\sigma_{j}\sigma_{k}).
		\end{displaymath}
		Using the commutation relations (\ref{commutations}), $\Pi(\delta\omega)$ reduces to 
		\begin{displaymath}
		(\sum_{i}\sum_{j=1}^{3}\partial_{j}(a_{i})\partial_{j}(b_{i})\ot\mathbb{I}
		-\sum_{i}a_{i}\partial_{1}(b_{i})\ot\sigma_{2}\sigma_{3}+\sum_{i}a_{i}\partial_{2}(b_{i})\ot\sigma_{1}\sigma_{3}++\sum_{i}a_{i}\partial_{3}(b_{i})\ot\sigma_{1}\sigma_{2}).
		\end{displaymath}
		But $\sum_{i}a_{i}\partial_{j}(b_{i})=0$ for $j=1,2,3$. So  
		\begin{displaymath}
		\Pi(\delta J_{0}^{1})\subset \clo_{3}\ot\mathbb{I}.
		\end{displaymath}
		We shall prove that the above inclusion is actually an equality. But this is straightforeward. Pick $\omega=S_{1}^{\ast}\delta S_{1}$. Then $S_{1}^{\ast}\partial_{j}(S_{1})=0$ for all $j$ so that $\Pi(\omega)=0$. But $\Pi(\delta\omega)=\sum_{j=1}^{3}\partial_{j}(S_{1}^{\ast})\partial_{j}(S_{1})\ot\mathbb{I}=2\ot\mathbb{I}$. Therefore,
		\begin{displaymath}
		\Pi(\delta J_{0}^{1})=\clo_{3}\ot\mathbb{I}.
		\end{displaymath}
		Hence we conclude that
		\begin{displaymath}
		\Omega_{\cld}^{2}(\clo_{3})=\clo_{3}\oplus\clo_{3}\oplus\clo_{3}.
		\end{displaymath}\qed.\vspace{0.2 in}\\
		From the proof of the Proposition \ref{1forms}, it follows that $\{1\ot\sigma_{i}:i=1,2,3\}$ is an $\clo_{3}$-basis for $\Omega^{1}_{\cld}(\clo_{3})$. We denote them by $\{e_{i}\}_{i=1,2,3}$. It is clear from the module structure that $e_{i}\in \mathcal{Z}(\Omega^{1}_{\cld}(\clo_{3}))$. Similarly $\{1\ot\sigma_{ij}\}_{i<j}$ is an $\clo_{3}$-basis for $\Omega^{2}_{\cld}(\clo_{3})$. We denote the basis elements by $\{e_{ij}\}_{i<j}$. With these notations we have the following
		\blmma
		\label{multiplication}
		The multiplication map $m:\Omega^{1}_{\cld}(\clo_{3})\ot_{\clo_{3}}\Omega^{1}_{\cld}(\clo_{3})\raro \Omega^{2}_{\cld}(\clo_{3})$ is surjective and is given by
		\begin{eqnarray}
		\label{mult}
		m(\sum_{i=1}^{3}e_{i}a_{i},\sum_{i=1}^{3}e_{i}b_{i})=\sum_{i<j}e_{ij}(a_{i}b_{j}-a_{j}b_{i}).
		\end{eqnarray}
		\elmma
		{\it Proof}:\\
		The surjectiveness of $m$ follows immediately once we prove (\ref{mult}). But (\ref{mult}) follows from  
		\begin{displaymath}
		m(e_{i}\ot e_{i})=0, m(e_{i}\ot e_{j})=e_{i}e_{j} \ for \ i<j, m(e_{i}\ot e_{j})=-e_{j}e_{i} \ for \ i>j.
		\end{displaymath}\qed\\
		From now on we shall denote the $\clo_{3}$-bimodule $\Omega^{1}_{\cld}(\clo_{3})$ by $\cle$.
		\bthm
		\label{existsunique}
		For the Connes' differential calculus coming from the $C^{\ast}$-dynamical system $(\clo_{3},\alpha, SO(3))$, a unique Levi-Civita connection exists for any pseudo-Riemannian metric. 
		\ethm
		{\it Proof}:\\
		We shall show that the Connes' differential calculus satisfies the assumptions {\bf I-IV} of subsection \ref{prel} and hence by Theorem \ref{gen}, the conclusion of the above theorem will follow. So let us check the assumptions one by one.\\
		{\bf Assumption I}: $\cle$ is a free module of rank $3$ with the basis elements $e_{1},e_{2},e_{3}\in\clz(\cle)$. Hence it is finitely generated and projective as both right and left module. It is easy to see that the map $u^{\cle}:\mathcal{Z}(\cle)\ot_{\clz(\clo_{3})}\clo_{3}\raro \cle$ given by $u^{\cle}(e_{i}\ot a_{i})=\sum e_{i}a_{i}$ is an isomorphism with the inverse $u^{\cle^{-1}}(\sum_{i=1}^{3}e_{i}a_{i})=\sum_{i=1}^{3}e_{i}\ot a_{i}$.\\
		{\bf Assumption II-III}: Recall the multiplication map $m$. From Lemma \ref{multiplication}, we have the following short exact sequence:
		
		\begin{displaymath}
		0\longrightarrow {\rm ker}(m)\longrightarrow \cle\ot_{\clo_{3}}\cle\stackrel{m}{\longrightarrow}\Omega^{2}_{\cld}(\clo_{3})\longrightarrow 0
		\end{displaymath}
		$\Omega^{2}_{\cld}(\clo_{3})$ being a free module, the above short exact sequence splits and hence we have
		\begin{displaymath}
		\cle\ot_{\clo_{3}}\cle={\rm ker}(m)\oplus \clf,
		\end{displaymath}
		where $\clf\cong \Omega^{2}_{\cld}(\clo_{3})$. As in \cite{levi}, This splitting implies that the map $P_{sym}$ on the elements $e_{i}\ot e_{j}$ is given by $P_{sym}(e_{i}\ot e_{j})=\frac{1}{2}(e_{i}\ot e_{j}+e_{j}\ot e_{i})$ and consequently $\sigma(e_{i}\ot e_{j})=e_{j}\ot e_{i}$ for all $i,j=1,2,3$. Thus Connes' differential calculus satisfies {\bf Assumptions II-III}.\\
		{\bf Assumption IV}: Consider the metric $g:\cle\ot_{\clo_{3}}\cle\raro\clo_{3}$ defined by 
		\begin{eqnarray}\label{metric}g(e_{i}\ot e_{j}a)=\delta_{ij}a, a\in\clo_{3}.\end{eqnarray}
		 Then it is easy to see that $g$ is a pseudo Riemannian bilinear non degenerate metric such that $g\circ\sigma=g$. Hence Assumption IV is satisfied. \qed.\\
		Henceforth we call the metric $g$ the {\bf canonical} metric.
		\section{Scalar curvature for the canonical metric}
		In this section we shall compute the scalar curvature of $\clo_{3}$ for the canonical metric $g$ given by (\ref{metric}).
		\subsection{Computation of Christoffel symbols}
		Let us recall from Theorem \ref{gen} that the Levi-Civita connection for any pseudo-Riemannian metric $g$ is given by $\nabla=\nabla_{0}+L$, where $\nabla_{0}$ is any torsionless connection (see definition \ref{torsionless}) and $L=\Phi_{g}^{-1}(dg-\Pi_{g}(\nabla_{0}))$. We prove the following Lemma which will enable us to choose a torsionless connection.
		\blmma
		\label{Tor1}
		\begin{eqnarray}
		&& de_{1}=m(e_{2}\ot e_{3})\\
		&& de_{2}=m(e_{3}\ot e_{1})\\
		&& de_{3}=m(e_{1}\ot e_{2})
		\end{eqnarray}
		\elmma
		{\it Proof}:\\
		Since $e_{1}=S_{2}^{\ast}d(S_{3})$, we have 
		\begin{eqnarray*}
		de_{1}&=&\overline{\Pi(\delta S_{2}^{\ast}\delta S_{3})}\\
		&=& (-S_{3}^{\ast}\ot\sigma_{1}-S_{1}^{\ast}\ot\sigma_{3})(S_{2}\ot\sigma_{1}+S_{1}\ot\sigma_{2})\\
		&=& 1\ot\sigma_{2}\sigma_{3}\\
		&=& m(e_{2}\ot e_{3})
		\end{eqnarray*}
		Similarly one can prove the other equalities.\qed
		\bppsn
		A choice of torsionless connection $\nabla_{0}$ is given on the $\clo_{3}$-basis of the free module $\Omega^{1}(\clo_{3})$ by
		\begin{eqnarray}
		\label{Tor2}
		\nabla_{0}(e_{1})=(e_{3}\ot e_{2}), \nabla_{0}(e_{2})=(e_{1}\ot e_{3}), \nabla_{0}(e_{3})=(e_{2}\ot e_{1}).
		\end{eqnarray}
		\eppsn
		{\it Proof}:\\
		Consequence of the Lemma \ref{Tor1} and the fact that $T_{\nabla_{0}}$ is right $\clo_{3}$-linear.\qed
		\blmma
		\label{2}
		For the canonical metric $g$ and the connection $\nabla_{0}$, let $\Pi_{g}(\nabla_{0})(e_{i}\ot e_{j})=-\sum_{m=1}^{3}T^{m}_{ij}e_{m}$. Then
		\begin{eqnarray}\
		\label{Chris}
		T^{m}_{ij}&=&-1, i,j,m \ are \ distinct \nonumber\\
		&=& 0, \ otherwise
		\end{eqnarray}
		\elmma
		{\it Proof}:\\
		Using equation (\ref{Tor2}),
		\begin{eqnarray*}
		\Pi_{g}(\nabla_{0})(e_{1}\ot e_{2})&=& (g\ot {\rm id})\sigma_{23}(\nabla_{0}(e_{1})\ot e_{2}+\nabla_{0}(e_{2})\ot e_{1})\\
		&=& (g\ot {\rm id})\sigma_{23}(e_{3}\ot e_{2}\ot e_{2}+ e_{1}\ot e_{3}\ot e_{1})\\
		&=& e_{3}
		\end{eqnarray*}
		With similar computation one can show that $\Pi_{g}(\nabla_{0})(e_{2}\ot e_{3})=e_{1}$ and $\Pi_{g}(\nabla_{0})(e_{1}\ot e_{3})=e_{2}$. Using the easy to see fact that $T^{m}_{ij}=T^{m}_{ji}$ for all $i,j,m$, we have for distinct $i,j,m$, $T^{m}_{ij}=-1$ for all $i,j,m=1,2,3$ and $T^{m}_{ij}=0$ for $i\neq j$ whenever either $m=i$ or $m=j$. So we are left with proving that $T^{m}_{ii}=0$ for all $i,m$. To that end
		\begin{eqnarray*}
		\Pi_{g}(\nabla_{0})(e_{1}\ot e_{1})&=& (g\ot{\rm id})\sigma_{23} 2(\nabla_{0}(e_{1})\ot e_{1})\\
		&=& 2(g\ot{\rm id})(e_{3}\ot e_{2}\ot e_{1})\\
		&=& 0.
		\end{eqnarray*}
		Similarly we can prove that $\Pi_{g}(\nabla_{0})(e_{i}\ot e_{i})=0$ for $i=2,3$, finishing the proof of the lemma.\qed
		\begin{thm}
		For the Levi-Civita connection $\nabla$, if we write $\nabla(e_{i})=\sum_{j,k=1}^{3}e_{j}\ot e_{k}\Gamma_{jk}^{i}$, for $i=1,2,3$, then we have the following:
		\begin{eqnarray}
		\label{Christoffel}
		\Gamma_{32}^{1}=\Gamma_{13}^{2}=\Gamma_{21}^{3}=.5, \Gamma_{23}^{1}=\Gamma_{31}^{2}=\Gamma_{12}^{3}=-.5.
		\end{eqnarray}
		Rest of the Christoffel symbols are zero.
		\end{thm}
		{\it Proof}:\\
		$\nabla-\nabla_{0}=L$ for some $L\in {\rm Hom}(\cle,\cle\ot^{\rm sym}_{\clo_{3}}\cle)$ (see Theorem 2.14 of \cite{levi}). If we write $L(e_{j})=\sum_{i,m=1}^{3}e_{i}\ot e_{m}L^{j}_{im}$ for $j=1,2,3$, then adapting the same arguments of Theorem 4.24 of \cite{levi},
		\begin{displaymath}
		L^{j}_{im}=\frac{1}{2}(T^{m}_{ij}+T^{i}_{jm}-T^{j}_{im}).
		\end{displaymath}
		Hence by equation (\ref{Chris}),
		\begin{eqnarray}
		\label{Chris2}
		L^{i}_{jm}&=&-0.5, \ i,j,m \ are \ distinct \nonumber\\
		&=& 0, \ otherwise.
		\end{eqnarray}
		
		\begin{eqnarray*}
		\nabla(e_{1})&=& \nabla_{0}(e_{1})+L(e_{1})\\
		&=& e_{3}\ot e_{2}-0.5 e_{3}\ot e_{2}-0.5 e_{2}\ot e_{3} \ (by \ (\ref{Chris2}))\\
		&=& 0.5 e_{3}\ot e_{2}-0.5 e_{2}\ot e_{3}.
		\end{eqnarray*}
		Hence $\Gamma^{1}_{23}=-.5, \Gamma^{1}_{32}=.5, \Gamma^{1}_{11}=\Gamma^{1}_{22}=\Gamma^{1}_{12}=\Gamma^{1}_{21}=0$. With similar computations, we can prove the rest of the equalities of (\ref{Christoffel}).\qed 
		\subsection{Computation of scalar curvature}
		Scalar curvature is obtained by contracting the Ricci tensor by the metric. Recall the maps $\rho$ and ${\rm ev}$ from Subsection \ref{prel}. The Ricci tensor ${\rm Ric}\in\cle\ot_{\clo_{3}} \cle$ is given by
		\begin{displaymath}
		{\rm Ric}=({\rm id}\ot {\rm ev}\circ\rho)(\Theta),
		\end{displaymath}
		where $\Theta$ is the curvature operator as in \cite{levi}.  
		\blmma
		$R(\nabla)(e_{i})=\frac{1}{8}\sum_{k\neq i}(e_{k}\ot e_{k}\ot e_{i}-e_{k}\ot e_{i} \ot e_{k})$, for $i=1,2,3$.
		\elmma
		{\it Proof}:\\
		Recall the map $H$ from (\ref{H}). By definition of $R(\nabla)$, we have
		\begin{eqnarray*}
			R(\nabla)(e_{1})&=& H\circ\nabla(e_{1})\\
			&=& \frac{1}{2}H(e_{3}\ot e_{2}-e_{2}\ot e_{3})\\
			&=& \frac{1}{2}[(1-P_{\rm sym})_{23}(\nabla(e_{3})\ot e_{2}-\nabla(e_{2})\ot e_{3})+e_{3}\ot Q^{-1}(de_{2})-e_{2}\ot Q^{-1}(de_{3})] \ (by \ (\ref{H}))
		\end{eqnarray*}
	We calculate each terms individually.
	\begin{eqnarray*}
		(1-P_{\rm sym})_{23}(\nabla(e_{3})\ot e_{2})&=&\frac{1}{2}(1-P_{\rm sym})_{23}(e_{2}\ot e_{1}\ot e_{2}-e_{1}\ot e_{2}\ot e_{2})\\
		&=& \frac{1}{4}(e_{2}\ot e_{1}\ot e_{2}-e_{2}\ot e_{2}\ot e_{1})
		\end{eqnarray*}
	Similarly $(1-P_{\rm sym})_{23}(\nabla(e_{2})\ot e_{3})=\frac{1}{4}(e_{3}\ot e_{3}\ot e_{1}-e_{3}\ot e_{1}\ot e_{3})$. Also using the fact that $Q^{-1}(e_{i}\ot e_{j})=\frac{1}{2}(e_{i}\ot e_{j}-e_{j}\ot e_{i})$ for all $i,j$, it is easy to see that
	\begin{displaymath}
		e_{3}\ot Q^{-1}(de_{2})=\frac{1}{2}(e_{3}\ot e_{3}\ot e_{1}-e_{3}\ot e_{1}\ot e_{3}), e_{2}\ot Q^{-1}(de_{3})=\frac{1}{2}(e_{2}\ot e_{1}\ot e_{2}-e_{2}\ot e_{2}\ot e_{1}).
		\end{displaymath}
		Hence 
		\begin{eqnarray*}
			R(\nabla)(e_{1})&=&\frac{1}{8}(e_{2}\ot e_{1}\ot e_{2}-e_{2}\ot e_{2}\ot e_{1}-e_{3}\ot e_{3}\ot e_{1}+e_{3}\ot e_{1}\ot e_{3}\\
			&& +2 e_{3}\ot e_{3}\ot e_{1}-2 e_{3}\ot e_{1}\ot e_{3}-2e_{2}\ot e_{1}\ot e_{2}+2e_{2}\ot e_{2}\ot e_{1})\\
			&=& \frac{1}{8}(e_{2}\ot e_{2}\ot e_{1}+ e_{3}\ot e_{3}\ot e_{1}-e_{2}\ot e_{1}\ot e_{2}-e_{3}\ot e_{1}\ot e_{3}).
			\end{eqnarray*}
		With similar computations, we have
		\begin{eqnarray*}
			&& R(\nabla)(e_{2})=\frac{1}{8}(e_{1}\ot e_{1}\ot e_{2}+ e_{3}\ot e_{3}\ot e_{2}-e_{1}\ot e_{2}\ot e_{1}-e_{3}\ot e_{2}\ot e_{3}),\\
			&&  R(\nabla)(e_{3})=\frac{1}{8}(e_{1}\ot e_{1}\ot e_{3}+ e_{2}\ot e_{2}\ot e_{3}-e_{1}\ot e_{3}\ot e_{1}-e_{2}\ot e_{3}\ot e_{2}).
			\end{eqnarray*}
		Combining these, we get the result of the lemma.\qed
		\blmma
		\label{Ricci}
		\begin{eqnarray*}
			{\rm Ric}(e_{i},e_{i})&=&-.25, \ for \ i=1,2,3\\
			&=& 0, \ i\neq j,
			\end{eqnarray*}
		i.e. ${\rm Ric}=\sum_{j=1}^{3}(-.25) e_{j}\ot e_{j}$.
		\elmma
		{\it Proof}:\\
		We begin by computing the curvature operator $\Theta$. To that end, recall the map $\zeta_{\cle,\clf}$ from (\ref{zeta}). Denoting the element of $\cle^{\ast}$ which takes the value $1$ at $e_{i}$ and $0$ on the other basis elements by $e_{i}^{\ast}$, we have
		\begin{eqnarray*}
			\Theta &=& (\sigma_{23}\ot {\rm id})\zeta_{\cle,\cle\ot_{\clo_{3}}\cle\ot_{\clo_{3}}\ot\cle} ^{-1}R(\nabla)\\
			&=& (\sigma_{23}\ot {\rm id})\frac{1}{8}\sum_{i,j:i\neq j}(e_{j}\ot e_{j}\ot e_{i}\ot e_{i}^{\ast}-e_{j}\ot e_{i}\ot e_{j}\ot e_{i}^{\ast})\\
			&=& \frac{1}{8}\sum_{i,j:i\neq j}(e_{j}\ot e_{i}\ot e_{j}\ot e_{i}^{\ast}-e_{j}\ot e_{j}\ot e_{i}\ot e_{i}^{\ast})
			\end{eqnarray*} 
		Now recall the map $\rho:\cle\ot\cle^{\ast}\raro\cle^{\ast}\ot \cle$. It is easy to see that $\rho(e_{i}\ot e_{j}^{\ast})=e_{j}^{\ast}\ot e_{i}$ for all $i,j=1,2,3$. Using this, we get
		\begin{eqnarray*}
			{\rm Ric}&=&({\rm id}_{\cle\ot\cle}\ot {\rm ev}\circ\rho)(\Theta)\\
			&=& \frac{1}{8}\sum_{i,j:i\neq j}({\rm id}_{\cle\ot\cle}\ot {\rm ev}\circ\rho)(e_{j}\ot e_{i}\ot e_{j}\ot e_{i}^{\ast}-e_{j}\ot e_{j}\ot e_{i}\ot e_{i}^{\ast})\\
			&=& \frac{1}{8}\sum_{i,j:i\neq j}(e_{j}\ot e_{i}{\rm ev}(e_{i}^{\ast}\ot e_{j})-e_{j}\ot e_{j}{\rm ev}(e_{i}^{\ast}\ot e_{i})).
			\end{eqnarray*}
		Since, by definition, ${\rm ev}(e_{i}^{\ast}\ot e_{j})=\delta_{ij}$, we have
		\begin{displaymath}
		{\rm Ric}=\sum_{j=1}^{3}(-.25) e_{j}\ot e_{j}.
		\end{displaymath}
		This proves the Lemma.\qed
		\bthm
		\label{Scal}
		Cuntz algebra with three generators with the differential calculus coming from the action of $SO(3)$ has a constant negative scalar curvature $-0.75$ for the canonical metric.
		\ethm{\it Proof}:\\
		By the definition of the scalar curvature (see after definition 4.4 in \cite{levi}),
		\begin{eqnarray*}
			{\rm Scal}&=&{\rm ev}\circ (V_{g}\ot {\rm id})({\rm Ric})\\
			&=& (-.25)\sum_{j=1}^{3}g(e_{j}\ot e_{j})\\
			&=&-.75
			\end{eqnarray*}\qed.
			\brmrk
			We have restricted our attention to Cuntz algebra with three generators purely because the computation for the general case does not seem to be that much more enlightening. \ermrk
			{\bf Acknowledgement}: The author acknowledges support from Department of Science and Technology, India (DST/INSPIRE/04/2016/002469). He would like to thank Jyotishman Bhowmick for some fruitful discussions and both K. B. Sinha as well as Debashish Goswami for their words of encouragement.
		
	\end{document}